\documentclass{article}

\usepackage{arxiv}

\usepackage[utf8]{inputenc}   % Allow UTF-8 input
\usepackage[T1]{fontenc}      % Use T1 fonts

\usepackage{amsmath, amssymb, amsthm, amsfonts}
\usepackage{mathtools}        % Better math (includes \DeclarePairedDelimiter)
\usepackage{mathrsfs}         % Script fonts

\DeclarePairedDelimiter\abs{\lvert}{\rvert}   % Absolute value

\usepackage{graphicx}
\graphicspath{{./images/}}    % Set path for images

\usepackage{hyperref}         % Hyperlinks
\usepackage{url}              % URL formatting
\usepackage{booktabs}         % Nicer tables
\usepackage{nicefrac}         % For nice ½ symbols
\usepackage{microtype}        % Better typography
\usepackage{lipsum}           % Placeholder text

\newtheorem{theorem}{Theorem}[section]
\newtheorem{definition}[theorem]{Definition}
\newtheorem{lemma}[theorem]{Lemma}
\newtheorem{proposition}[theorem]{Proposition}
\newtheorem{corollary}[theorem]{Corollary}

\DeclarePairedDelimiter\norm{\lVert}{\rVert}
\usepackage{xcolor}  % already included via mathtools dependencies
\usepackage{graphicx}
\graphicspath{{./images/}}    % Set path for images
\usepackage{float}            % For [H] placement
\usepackage{caption}          % Better control over captions
\usepackage{subcaption}
\DeclareMathOperator{\supp}{supp}

\newcommand{\note}[1]{\textcolor{black}{\textbf{[Note :} #1\textbf{]}}}

\title{Quasi-compactness of Frobenius-Perron operator for piecewise expanding $C^{1+\varepsilon}$ maps of an interval}

\author{
 Aparna Rajput \\
  Department of Mathematiocs and Statistics\\
  Concordia University\\
  Montreal, H3G 1M8, QC, Canada\\
  \texttt{a\_ajpu@live.concordia.ca} \\
   \And
 Pawe\l \ G\'ora \\
 Department of Mathematiocs and Statistics\\
  Concordia University\\
  Montreal, H3G 1M8, QC, Canada\\
  \texttt{pawel.gora@concordia.ca} \\
}

\begin{document}
\maketitle
\begin{abstract}
In this paper, we establish  Lasota-Yorke inequality for the Frobenius-Perron Operator of a piecewise expanding $C^{1+\varepsilon}$ map of an interval. By adapting this inequality to satisfy the assumptions of the Ionescu-Tulcea and Marinescu ergodic theorem, we demonstrate the existence of an absolutely continuous invariant measure (ACIM) for the map. Furthermore, we prove the quasi-compactness of the Frobenius-Perron operator induced by the map. Additionally, we explore significant properties of the system, including weak mixing and exponential decay of correlations.
\end{abstract}

% keywords can be removed
\keywords{Absolutely continuous invariant measures \and  Exactness of a dynamical system \and  Piecewise expanding map}


\bigskip
\noindent \textbf{2020 Mathematics Subject Classification.} 37A05; 37E05.

\section{Introduction}
In this paper, we study the piecewise expanding, piecewise \( C^{1+\varepsilon} \) maps \( \tau \) of the unit interval. Our goal is to prove the quasi-compactness of the Frobenius-Perron operator \( P_\tau \) on the spaces \( BV_{1,\frac 1p} \), where \( p \geq 1 \) and \( 1/p = \varepsilon \). In Section 4 we also prove some results on the behaviour of $P_\tau$ on the spaces $BV_{t,\frac 1p},
1\le t \le p$, which, in the future, may be useful for showing the quasi-compactness on these spaces.

These spaces of functions of generalized bounded variation were introduced by Keller \cite{keller1985}. He proved the quasi-compactness of \( P_\tau \) on \( BV_{1,\frac 1p} \) under the weaker assumption that \( \frac {1}{|\tau'|} \in BV_{1,\frac 1p} \). This assumption would not allow us to generalize the result to  $BV_{t,\frac 1p}$, $t>1$, as we need the
boundedness of $\tau'$.\\
We work within the space \( L^1 \) of all integrable functions defined on the unit interval \( I \) with respect to the Lebesgue measure \( m \).
A measurable transformation \( \tau: I \rightarrow I \) is called nonsingular if, for any measurable set \( A \subset I\), the condition \( m(A) = 0 \) implies that \( m(\tau^{-1}(A)) = 0 \) as well. The piecewise $C^{1+\epsilon}$ expanding maps are nonsingular.

For $\tau$ piecewise monotonic, the Frobenious-Perron Operator 
$\displaystyle P_\tau:L^1 \rightarrow L^1$ is defined by,
\begin{equation*}
    P_\tau f(x)=\sum_{y\in \tau^{-1}(x)} f(y) \cdot g(y) =\sum_{i=1}^N \frac{f(\tau^{-1}(x)}{\abs{\tau'(\tau^{-1}(x))}} \cdot \chi_{\tau(I_i)}(x),
\end{equation*}
 where $g(y)=\frac{1}{\abs{\tau'(y)}}$ and $\sup_{I_i} \abs{g}<\frac{1}{s_i}<1$.
We know that $P_\tau$ is a linear and continuous operator and has the following properties:\\
(1) $P_\tau$ is positive, It means $f\geq 0 \implies P_\tau f \geq 0$;\\
(2) $P_\tau$ preserves integrals, i.e. $\displaystyle \int_I P_\tau f dm =\int_I f dm$ for $f\in L_1$;\\
(3) $P_\tau$ satisfies the composition property, $\displaystyle P_{\tau^n} =P^n_\tau$ where $n$ is the nth iterate of $\tau$;\\
(4) $P_\tau f= f$ if and only if the measure $d\mu=f dm$ is $\tau$ invariant, i.e.  for each measurable set $A, \mu(\tau^{-1} (A))=\mu(A)$. The Perron-Frobenius operator \( P_\tau \) is a powerful tool for capturing the ergodic properties of the dynamical system \( (\tau, \mu) \), particularly in terms of :

\begin{itemize}
    \item Existence of an absolutely continuous invariant measure (ACIM),
    \item Weak mixing and decay of correlations,
     \item Quasi-compactness.
\end{itemize}
We reprove Keller's result under the assumption that $\tau \in C^{1+\varepsilon}$, 
with the estimate $\frac{1}{s^{1/p}} + \frac{2}{s} < 1$ on the minimal slope of $\tau$. 
Keller's estimate is $\frac{2}{s^{1/p}} < 1$ (Lemma~3.1 of \cite{keller1985}). We establish Lasota-Yorke inequality for the pair $(BV_{1,\frac{1}{p}}, L^1)$ in Theorem \ref{thm 3.2.2}.  Theorem \ref{thm 3.2.3} shows that the iterates of $P_\tau$ are uniformly bounded in $BV_{1,\frac{1}{p}}$ norm. The quasi-compactness of \( P_\tau \) has a number of interesting consequences. It demonstrates that the map \( \tau \) has a finite number of ergodic components, and on each component, it possesses an absolutely continuous invariant measure, see  Theorem \ref{thm 3.2.4}. Furthermore, for some iterate of \( \tau \), say \( \tau^n \), there are finitely many disjoint invariant domains on which \( \tau^n \) is exact. In Theorem \ref{Thm 3.2.7}, we establish that for the system \((\tau, \mu)\), which admits a unique absolutely continuous invariant measure \(\mu\), the correlation \(C_m(f, g, N)\) decays exponentially. Moreover, Theorem \ref{thm 3.2.8} shows exponential decay for more symmetric form of correlation computed with respect to the invariant measure \(\mu\).

There was a strong interest in $C^{1+\varepsilon}$ maps in the 1970s, when the chaos in the Lorenz
system was not proven yet. It came from the fact that a possibly important map is only
$C^{1+\varepsilon}$ smooth and not $C^2$. This map, suggested by Lorenz himself~\cite{lorenz2017},
is ``the next maximum of the $z$-map in a typical
trajectory of the Lorenz system,'' see Figure~(b). It is constructed as follows~\cite{fang1995,boyarsky2015}.
From a trajectory of the Lorenz system, we take consecutive points for which the $z$-coordinate attains
a local maximum (in the trajectory), see Figure~(a). The figure has been prepared using the Runge--Kutta (RKF45) approximation algorithm~\cite{maplesoft2017}.
The consecutive values of $z$ define a map that resembles a tent map
but is believed to be only $C^{1+\varepsilon}$. Proving that this map is chaotic (e.g., has an ACIM) was considered to be a strong indication of
chaos in the Lorenz system.

Since then, the chaotic behaviour of the Lorenz system has been rigorously established; 
see, for example,~\cite{afraimovich1977,guckenheimer1979,tucker1999}.
\begin{figure}[H]
  \centering
  \begin{subfigure}[b]{0.50\textwidth}
    \includegraphics[width=\textwidth]{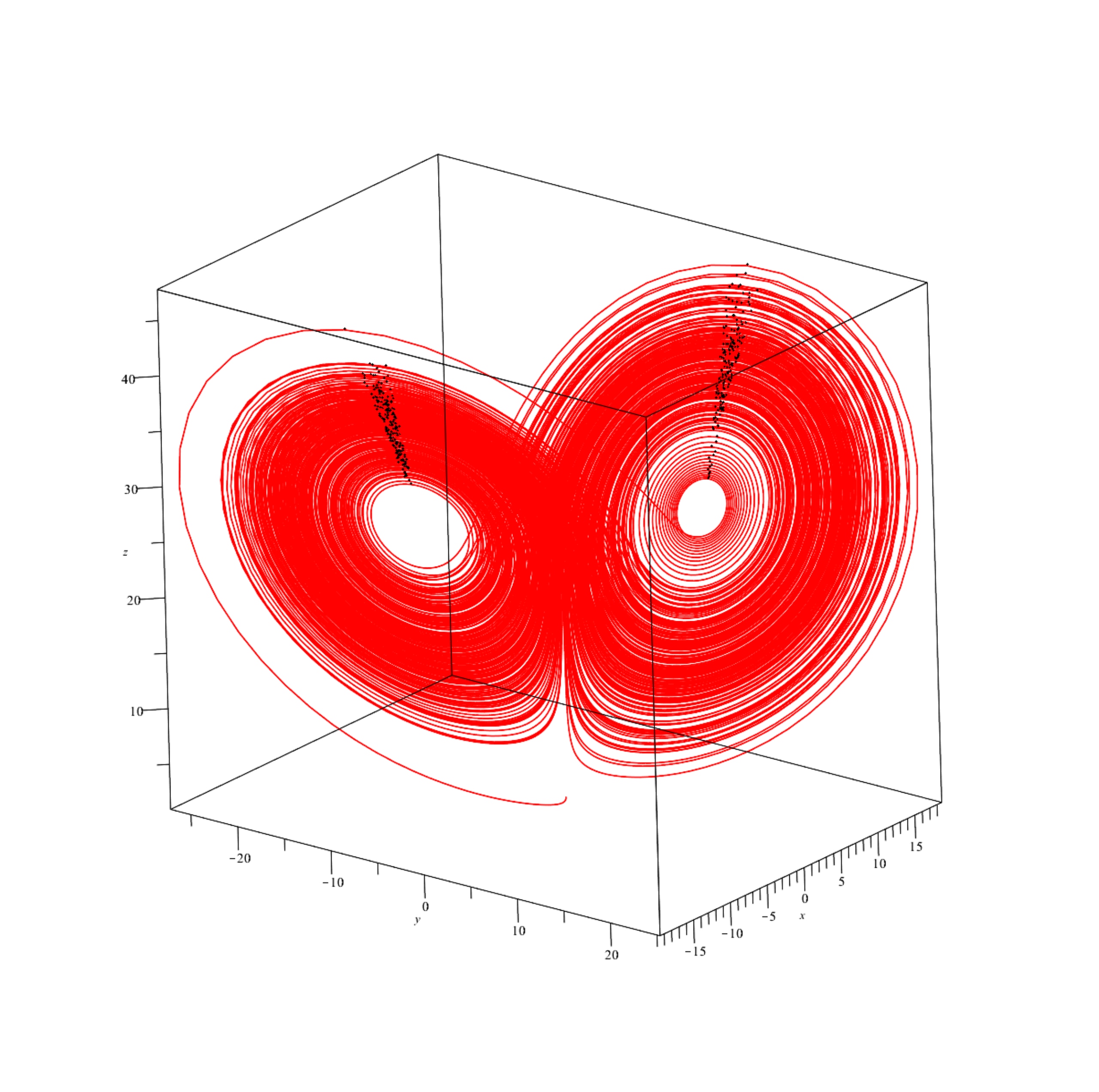}
    \caption{A plot of the orbit starting at $(0.1, -0.3, 1.7)$ for the Lorenz system. 
The set of local maxima points is shown in black.}
  \end{subfigure}
  \hfill
  \begin{subfigure}[b]{0.48\textwidth}
    \includegraphics[width=\textwidth]{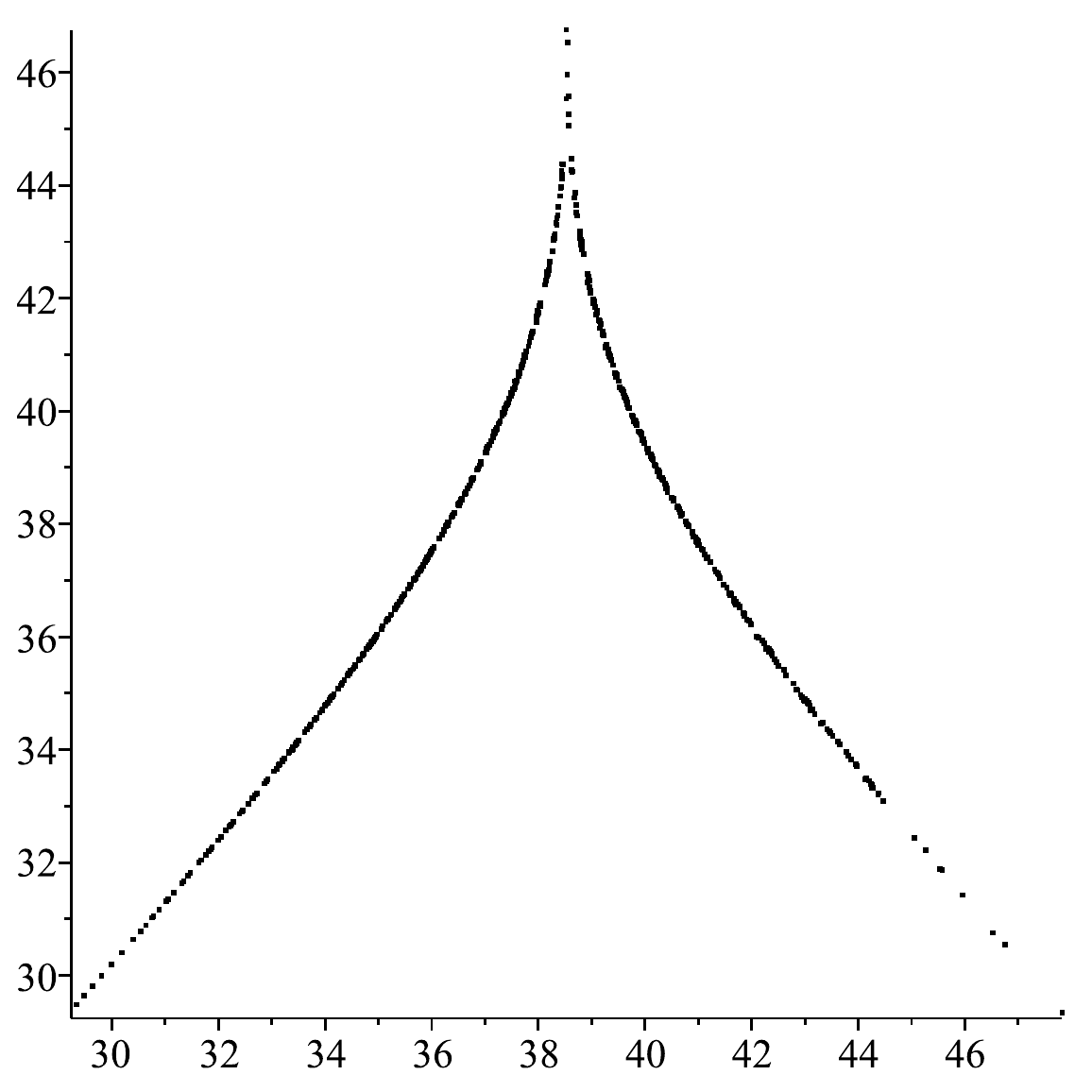}
    \caption{The successive maxima map plotted using the points from Figure~(a).}
  \end{subfigure}
\end{figure}
\section[Generalized bounded variation and applications to piecewise monotonic transformations]{Generalized bounded variation and applications to \mbox{piecewise} monotonic transformations}\label{sec 1.7}
Keller \cite{keller1985} used the concept of generalized bounded variation in application to piecewise monotonic transformations, when the inverse of the derivative is Hölder-continuous or, more generally, of bounded \( p \)-variation on the unit interval $I=[0,1]$. He investigated the quasi-compactness of the Perron-Frobenius operator and the existence of absolutely continuous invariant measure.
\begin{definition}\label{def1.7.1}
    (\text{osc}illation of a function)
    For an arbitrary function $f: I \rightarrow \mathbb{C} $ and $r>0$, we define \text{osc}$(f,r,.): I \rightarrow [0,\infty)$ by,
    \[
\text{\text{osc}}(f, r, x) = \sup_{y_1, y_2 \in S_r(x)} |f(y_1) - f(y_2)|,
\]
where \( S_r(x) = (x - r, x + r) \cap I \).
\end{definition}
The function $\text{\text{osc}}(f, r, x)$ measures how much the function $f$ oscillates within the neighborhood  $S_r (x)$, essentially capturing the largest difference in the values of $f$ within this ball.
\begin{proposition}\label{prop 1.7.2}
  The function \( \mathrm{\text{osc}}(f, r, .) \)  is lower semi-continuous.
\end{proposition}
\begin{proof}
    Fix \( r > 0 \) and \( x_0 \in I \). Given \( \varepsilon > 0 \), there exist \( y_1, y_2 \in S_r(x_0) \) such that 
\[
|f(y_1) - f(y_2)| > \text{\text{osc}}(f, r, x_0) - \varepsilon.
\]
Let \( 0 < \delta < r - \max\{|y_1 - x_0|, |y_2 - x_0|\} \). With this, if \( x \in I \) and \( |x - x_0| < \delta \), then \( y_1, y_2 \in S_r(x) \), and consequently,
\[
\text{\text{osc}}(f, r, x) \geq |f(y_1) - f(y_2)| > \text{\text{osc}}(f, r, x_0) - \varepsilon.
\]
In other words,
\[
\text{\text{osc}}(f, r, x) > \text{\text{osc}}(f, r, x_0) - \varepsilon \quad \forall x \in S_\delta(x_0) \cap I,
\]
so \( \text{\text{osc}}(f, r, \cdot) \) is lower semi-continuous at \( x_0 \).

\end{proof}
Since lower semi-continuous functions are measurable, we conclude that \( \mathrm{\text{osc}}(f, r, .) \) is also a measurable function, and we define for $1\leq t\leq \infty $, $\text{\text{osc}}_t(f,r)=\|\text{\text{osc}}(f,r,.)\|_t$, where we allow the $t$ norm to take the value $+\infty$.
$\text{\text{osc}}_t(f,r)$ is an increasing function (in the variable $r$) from $(0, A]$ to $[0, \infty]$, for any positive constant $A$.\\

 \begin{definition}\label{def 1.7.3}
     For $t,p \geq 1$, we define:
     \begin{itemize}
     \item[(1)] For $f: I \to \mathbb{C}$, set $\text{var}_{t,\frac{1}{p}}$ variation of $f$ as 
    \[
    \text{var}_{t, \frac{1}{p}}(f) = \sup_{0 < r\leq A} \frac{\text{osc}_t(f,r)}{r^\frac{1}{p}};
    \]
    \item[(2)] $BV_{t, \frac{1}{p}}$ is the space of $m$-equivalence classes of functions of bounded $\text{var}_{t,\frac{1}{p}}$;

    \item[(3)] For $f \in BV_{t, \frac{1}{p}}$, we define the norm:
    \[
    \|f\|_{t,\frac{1}{p}} = \text{var}_{t, \frac{1}{p}}(f) + \|f\|_t.
    \]
\end{itemize}
 \end{definition}
Observe that the definition depends on the constant $A$.\\
For the proof of the following two results, we refer the reader to \cite{keller1985} and \cite{hofbauer1982}.
\begin{proposition}
    $BV_{t,\frac{1}{p}}$ is a linear space, and $\|.\|_{t,\frac{1}{p}}$ is a norm on it for $1\leq t, p \leq \infty$.
\end{proposition}
\begin{proposition}
 Let $\{f_n\}$ be a sequence in $BV_{t,\frac{1}{p}}$ converging in $\|.\|_t$ norm to some function $f\in L^t$. Then 
    \[
        \textnormal{var}_{t,\frac{1}{p}}(f)\leq \liminf_{n\rightarrow \infty}{\textnormal{var}_{t,\frac{1}{p}}(f_n)} .
    \]
\end{proposition}
\begin{proposition}\label{prop 1.7.6}
If \( f \in BV_{t,\frac 1p} \), then \( f \) is bounded almost everywhere and has a version which is lower semi-continuous.
\end{proposition}
\begin{proof}
Let \( f \in BV_{t,\frac 1p} \). First, we prove that \( f \) is bounded a.e.. We
prove that it is bounded above a.e.. Assume that it does not hold.
Then, for any \( N \geq 1 \) there is a set of positive measure \( A_N \) such that
\( f(x) \geq N \) for \( x \in A_N \). The sets \( A_N \) are decreasing, i.e., \( A_{N+1} \subset A_N \),
for all \( N \). Since they all are the subsets of \([0, 1]\), their closures intersect
at least at one point, say \( y \). Then, \( \operatorname{osc}(f,r, x) \) is infinite for any \( r \) and
any \( x \in S(y, r/2) \). This implies that
\[
\int_I \operatorname{osc}(f,r, x)^t \, dm = \infty
\]
for all \( r \), and \( f \) cannot be in \( BV_{t,\frac 1 p} \).

Next, we prove  lower semi-continuity. We will show that \( f \) is continuous a.e..
Let \( B \) be the set of discontinuity points of \( f \), i.e., for any \( x \in B \) at least
one of the one-sided limits does not exist or is not equal to \( f(x) \). Let us assume that \( m(B) > 0 \). Let \( x \in B \). Then, there exist sequences \( x_n \to x \) and \( y_n \to x \), as \( n \to \infty \), with
\[
\left| \lim_{n \to \infty} f(x_n) - \lim_{n \to \infty} f(y_n) \right| \geq a_x > 0.
\]
 For one sequence we can choose a constant sequence $x_n=x$, $n=1,2,\dots,$ with $\lim_{n\rightarrow \infty} f(x_n)=f(x)$. For the second one, we choose any sequence $\left\{y_n\right\}_{n=1}^\infty$ with $\lim_{n\rightarrow \infty} y_n=x$ and $\lim_{n\rightarrow \infty} f(y_n)=a \neq f(x)$, which must exist since at least one of the one-sided limits does not exist or is not equal to $f(x)$. This implies that \( \operatorname{osc}(f, x,r) \) is greater than \( a_x \), for all \( r > 0 \). Thus, we can find an \( n > 1 \)
and a set of positive measure \( B_n \) such that \( \operatorname{osc}(f, x,r) > 1/n \) for \( x \in B_n \)
and \( r > 0 \). Then, for \( r > 0 \) we have
\[
\int_I \operatorname{osc}(f, x,r)^t \, dm \geq \frac{m(B_n)}{n^t}.
\]

Then,
\[
\sup_{0<r \leq A} \frac{1}{r^p} \int_I \operatorname{osc}(f, x,r)^t \, dm
\]
is infinite and \( f \) can not be in \( BV_{t, \frac 1 p} \).

Thus, the set \( B \) on which \( f \) is not continuous is of measure 0. If
we replace values of \( f \) for \( x \in B \) by \( \liminf_{y \to x} f(y) \), we obtain a
lower semi-continuous version \( \tilde{f} \) of \( f \).

Finally, we show that this modification does not damage the continuity
of \( \tilde{f} \) outside the set \( B \). Let \( x_0 \) be a continuity point of \( f \) and assume \( x_n \to x_0 \) as \( n \to \infty \). Let \( L_n =\tilde{f}(x_n)= \liminf_{y \to x_n} f(y) \). For any
\( n = 1, 2, \ldots \) we can find a point \( y_n \) such that \( |y_n - x_n| < 1/n \) and
\( |f(y_n) - L_n| < 1/n \). Then, \( y_n \to x_0 \) and since \( x_0 \) is a continuity
point of \( f \),
\[
f(x_0) = \lim_{n \to \infty} f(y_n) = \lim_{n \to \infty} L_n.
\]
This shows that \( \tilde{f}(x_0) = \lim_{n \to \infty} \tilde{f}(x_n) \), and \( \tilde{f} \) is continuous at all points where \( f \)
is continuous.
\end{proof}
We refer to a theorem from \cite{ionescu1950}.
\begin{theorem}\label{thm 1.7.7}
    For $1\leq  t,p < \infty$ we have :
    \begin{enumerate}
      \item[(1)] $E=\{f\in BV_{t,\frac{1}{p}} \mid \|f\|_{t,\frac{1}{p}} \leq c\}$ is a compact subset of $L^t \, $ for each $c>0$ ;
         \item[(2)] $(BV_{t,\frac{1}{p}}, \|.\|_{t,\frac{1}{p}})$ is a Banach space ;
     \item[(3)]$BV_{t,\frac{1}{p}}$ is dense in $L^t (1\leq t <\infty)$.
    \end{enumerate}
\end{theorem}
\begin{proof}
    We refer \cite{ionescu1950} for the proof.
\end{proof}
We will be using the following estimate given by \cite{keller1985}.
\begin{proposition}\label{prop 1.7.8}
   Let \( X, Y \) be subsets of \( I \), and let \( \tau: X \to Y \) be a bijection with $\abs{\tau'}\geq s$. Suppose \( f: X \to \mathbb{C} \) is a function. Then,
\[
\operatorname{osc}(f \circ \tau^{-1}, r, y) \leq \operatorname{osc}(f, \frac{r}{s}, \tau^{-1}y)
\]
for each \( y \in Y \).
\end{proposition}
We will use a version of Lemma 1.8, \cite{keller1985}.
 \begin{proposition}\label{prop 1.7.9}
Let $Y$ be a closed interval of length $kB$, $k=1,2,\ldots$. Then, for 
$f \in BV_{t,\frac 1p}$ we have
\[
\sup_{z \in Y} |f(z)| \;\le\; \frac{1}{B}\int_Y \operatorname{osc}(f,B,x)\,dm(x) 
\;+\; \frac{1}{B}\int_Y |f(x)|\,dm(x).
\]
\end{proposition}

\begin{proof}
Let the points $y_0 < y_1 < \cdots < y_k$ be the points dividing $Y$ into 
$k$ subintervals $Y_i = [y_{i-1},y_i]$ of equal length. 
Let $|f(z_0)| = \inf_{z\in Y} |f(z)|$. Since $f$ can be considered lower 
semi-continuous, the infimum is attained. 

Let $z\in Y$, $z_0<z$. We can find consecutive points 
$y_i,y_{i+1},\ldots,y_{i+t}$ such that 
$z_0,y_i\in Y_i$, $y_i,y_{i+1}\in Y_{i+1}, \ldots, y_{i+t}, z\in Y_{i+t}$, and
\[
|f(z)| \;\le\; |f(z_0)| + |f(z_0)-f(y_i)| + |f(y_i)-f(y_{i+1})| 
+ \cdots + |f(y_{i+t})-f(z)|.
\]
For any $x\in Y_i$ we have $|f(z_0)-f(y_i)| \le \operatorname{osc}(f,B,x)$,  
for any $x\in Y_{i+1}$ we have $|f(y_i)-f(y_{i+1})| \le \operatorname{osc}(f,B,x)$,  
and similarly, for any $x\in Y_{i+t}$ we have $|f(y_{i+t})-f(z)| \le \operatorname{osc}(f,B,x)$.  
Thus,
\[
|f(z)| \;\le\; \frac{1}{B}\int_Y \operatorname{osc}(f,B,x)\,dm(x) 
\;+\; \frac{1}{B}\int_Y |f(x)|\,dm(x).
\]
The proof for $z<z_0$ is similar.
\end{proof}
\begin{proposition}\label{Prop 1.7.10}
Under the assumptions of Proposition~\ref{prop 1.7.9} we have
\begin{equation}\label{eq:prop02}
\sup_{z_1,z_2\in Y} |f(z_1)-f(z_2)| 
\;\le\; \frac{1}{B}\int_Y \operatorname{osc}(f,B,x)\,dm(x).
\end{equation}
\end{proposition}
\begin{proof}
The proof is the same as for Proposition~\ref{prop 1.7.9}.
\end{proof}
\begin{proposition}\label{prop 1.7.10}
    Under the assumptions of Proposition \ref{prop 1.7.9}, if $J \subset Y$ is a subinterval of $Y$ and $t>1$ then 
    \begin{equation*}
        \left(\int_J \abs{f(z)}^t \, dm(z)\right)^\frac{1}{t} \leq \left(\frac{m(J)}{B}\right)^\frac{1}{t} \left[\left(\int_Y \left(\operatorname{osc}(f, B, x)\right)^t\, dm(x)\right)^\frac{1}{t} + \norm{f_{\vert Y}}_t\right] .
    \end{equation*}
\end{proposition}
\begin{proof}
    We estimate
    \begin{equation*}
    \begin{aligned}
       \left(\int_J \abs{f(z)}^t \, dm(z)\right)^\frac{1}{t} &\leq \left(\int_J \left(\frac{1}{B} \int_Y \operatorname{osc}(f, B, x) \, dm(x)  +\frac{1}{B} \int_Y |f| \, dm(x)\right)^t\, dm(z)\right)^\frac{1}{t} &\\
       &\leq \left(\int_J \left(\frac{1}{B} \int_Y \operatorname{osc}(f, B, x) \, dm(x) \right)^t\, dm(z)\right)^\frac{1}{t}+\left(\int_J\left(\frac{1}{B} \int_Y |f| \, dm(x)\right)^t\, dm(z)\right)^\frac{1}{t}.
    \end{aligned}
    \end{equation*}
    By Jensen's inequality we have, 
    \begin{equation*}
     \left(\frac{1}{B} \int_Y \operatorname{osc}(f, B, x) \, dm(x) \right)^t\leq    \frac{1}{B} \int_Y \left(\operatorname{osc}(f, B, x)\right)^t \, dm(x),
    \end{equation*}
    and
    \begin{equation*}
        \left(\frac{1}{B} \int_Y |f| \, dm(x)\right)^t \leq \frac{1}{B} \int_Y |f|^t \, dm(x).
    \end{equation*}
    Hence,
    \begin{equation*}
        \begin{aligned}
         \left(\int_J \abs{f(z)}^t \, dm(z)\right)^\frac{1}{t} &\leq   \left(\frac{m(J)}{B}\right)^\frac{1}{t} \left(\int_Y \left(\operatorname{osc}(f, B, x)\right)^t \, dm(x)\right)^\frac{1}{t}+\left(\frac{m(J)}{B}\right)^\frac{1}{t} \left(\int_Y |f|^t \, dm(x) \right)^\frac{1}{t}&\\
         &= \left(\frac{m(J)}{B}\right)^\frac{1}{t} \left[\left(\int_Y \left(\operatorname{osc}(f, B, x)\right)^t \, dm(x)\right)^\frac{1}{t} + \norm{f_{\vert Y}}_t\right] .
        \end{aligned}
    \end{equation*}
\end{proof}
\section{Piecewise expanding \texorpdfstring{$C^{1+\varepsilon}$}{C¹⁺ᵋ} maps of an interval}
\begin{definition}
    Assume there exists a finite partition $\mathcal{P}=\{I_i=(a_{i-1},a_i),i=1,2,3...q\}$ of $I=[0,1]$ and $\varepsilon>0$ such that the transformation $\tau:I\rightarrow I$ satisfies the following conditions:
\begin{enumerate}
    \item $\tau_{|_{I_i}}$ is a monotonic, $C^1$ function, which can be extended to a $C^1$ function $\tau_i$ on $\overline{I_i}$ for each $ i=1,2,...,q$;
    \item $\tau_i' (x)$ is $\varepsilon$-H\"older continuous for each  $i=1,2,...,q$, i.e., there exist a constant $M_i$ such that $|\tau_i' (x)-\tau_i' (y)|\leq M_i |x-y|^\varepsilon$, for all $x,y \in I_i$ and in particular if $\varepsilon=1$, then $\tau_i'$ is a Lipschitz function;
\item $|\tau_i' (x)|\geq s_i >1$ for all $x\in I_i$.
\end{enumerate}
We denote the class of all maps $\tau:I\rightarrow I$ satisfying $(1)-(3)$ by $\mathcal{T}_{1+\varepsilon}$, the class of piecewise expanding $C^{1+\varepsilon}$ maps on $I$.
\end{definition}
\note\textnormal{ If  $\varepsilon>1$, $\tau_i'(x)$ satisfies  H\"older's condition of order $\varepsilon>1$, $\tau_i(x)$ will be a constant function.
For $\varepsilon>1$, each $\tau_i'$ is a constant function and hence $\tau$ is piecewise linear.}\\

\textnormal{Now, we establish the Lasota-Yorke inequality, for $f\in BV_{1,\frac{1}{p}}$, for any $p\geq1$.}
\begin{lemma}\label{lem 3.2.1}
    Let $\tau \in \mathcal{T}_{1+\varepsilon}$ for some $p\geq1$ and $\varepsilon=1/p$. Then there exist $\alpha, \beta >0$ such that, for every $f\in BV_{1,\frac{1}{p}}$,
    \begin{equation}
       \textnormal{var}_{1,\frac{1}{p}}(P_\tau f) \leq \alpha \cdot \textnormal{var}_{1,\frac{1}{p}}(f) + \beta\cdot \norm{f}_1.
    \end{equation}
\end{lemma}
\begin{proof}
Recall from the definition of the Frobenious-Perron Operator  and Definition \ref{def1.7.1}  that, $S_r(x)=(x-r,x+r)\cap I$, and  $g(x)=\frac{1}{\abs{\tau_i'(x)}}$. We have
\begin{equation*}
\begin{aligned}
\text{osc}\,_1(P_\tau f,r)&=\int_I \text{osc}(P_\tau f,r,x) \,dm(x) &\\
   &= \displaystyle \int_I \sup_{y_1,y_2 \in S_r(x)}\abs{P_\tau f(y_1)-P_\tau f(y_2)} \,\, dm(x)&\\
&=\int_I\sup_{y_1,y_2 \in S_r(x)}\left|\sum_{i=1}^q \left((f\cdot g)(\tau_i^{-1}(y_1))\cdot\chi_{\tau(I_i)}(y_1)-(f\cdot g)(\tau_i^{-1}(y_2))\cdot \chi_{\tau(I_i)}(y_2)\right)\right|\,\, dm(x).
    \end{aligned}
\end{equation*}
For $i=1,2, \dots, q$ : if $\tau_i(a_{i-1})\notin \{0,1\}$ we define, $E_i^L=[\tau_i(a_{i-1})-r, \tau_i(a_{i-1})+r]\cap [0,1]$ and $F_i^L=E_i^L \cap\tau_i(I_i)$. If $\tau_{i}(a_i) \notin \{0,1\}$, then we define $E_i^R=[\tau_{i}(a_i)-r,\tau_{i}(a_i)+r]\cap [0,1]$ and $F_i^R=E_i^R \cap \tau_{i}(I_{i})$.
Note that $\tau_i^{-1}(E_i^L)=\tau_i^{-1}(F_i^L)$ and $\tau_i^{-1}(E_i^R)=\tau_i^{-1}(F_i^R)$. Let $$
F_{i}^{R\cup} := \tau(I_i)\cap B(E_i^R,r)
\quad\text{and}\quad
F_{i}^{L\cup} := \tau(I_i)\cap B(E_i^L,r),
$$
where, \(B(C, r)\) denotes the ball of radius $r$ around the set $C, i \in Q$. If, for example, $E_i^R=[c_1, c_2]$, $E_i^L=[c_3, c_4]$ and \(\tau_i\) is decreasing, then $F_{i}^{R\cup}=(c_1 - r,\, c_2]$ and $F_{i}^{L\cup}=[c_3,\, c_4 + r)$.\\
 Let $G_i=\tau(I_i)\setminus (E_i^R \cup E_i^L)$ and $h_i=(f\cdot g)\circ \tau_i^{-1}$, $i=1,\dots,q$. We have
\begin{equation*}
    \begin{aligned}
      \text{osc}\,_1(P_\tau f,r)&  \leq  \sum_{i=1}^q\int_{G_i}\sup_{y_1,y_2 \in S_r(x)}\left|h_i(y_1)-h_i(y_2)\right|\, dm(x)&\\
        &+ \sum_{i\in Q} \int_{E_i^R} \sup_{y_1, y_2\in F_{i}^{R\cup}} \left|h_i(y_1)-h_i(y_2) \right|\, \, dm(x)&\\
        &+\sum_{i\in Q}\int_{E_i^L} \sup_{y_1,y_2\in F_{i}^{L\cup}} \left| h_i(y_1)-h_i(y_2) \right|\,\, dm(x),
    \end{aligned}
\end{equation*}
where $Q$ is the set of indices $i$ such that $\tau_i$ is not onto $I$. We split this into two parts
\begin{equation}
    \Delta_1=\sum_{i=1}^q\int_{G_i} \sup_{y_1,y_2 \in S_r(x)} \left| h_i(y_1)-h_i(y_2)\right|\,\, dm(x).
\end{equation}
and
\begin{equation}
    \Delta_2= \sum_{i\in Q} \int_{E_i^R} \sup_{y_1, y_2\in F_{i}^{R\cup}} \left| h_i(y_1)-h_i(y_2) \right|\, \, dm(x)+ \sum_{i\in Q} \int_{E_i^L} \sup_{y_1,y_2\in F_{i}^{L\cup}} \left| h_i(y_1)-h_i(y_2) \right|\, \, dm(x).
\end{equation}
First, we estimate $\Delta_1$.
Adding and subtracting $f(\tau_i^{-1}(y_1))\cdot g(\tau_i^{-1}(y_2))$ we get,
\begin{equation}\label{equation 3.1.4}
\begin{aligned}
     \Delta_1 &=\sum_{i=1}^q \int_{G_i}\sup_{y_1,y_2 \in S_r(x)} \left|(f(\tau_i^{-1}(y_1))\bigl(g(\tau_i^{-1}(y_1))-g(\tau_i^{-1}(y_2))\right|&\\
   &+ g(\tau_i^{-1}(y_2))\bigl(f(\tau_i^{-1}(y_1))
    -f(\tau_i^{-1}(y_2))\bigl)\bigg| \,\, dm(x)&\\
    &\leq\sum_{i=1}^q \int_{G_i}\sup_{y_1,y_2 \in S_r(x)}  \left|f(\tau_i^{-1}(y_1))\bigl(g(\tau_i^{-1}(y_1))-g(\tau_i^{-1}(y_2))\bigl)\right|\,dm(x) &\\
          &+\sum_{i=1}^q \int_{G_i}\sup_{y_1,y_2 \in S_r(x)} \left| g(\tau_i^{-1}(y_2))\bigl(f(\tau_i^{-1}(y_1))-f(\tau_i^{-1}(y_2))\bigl)\right|\,dm(x).
\end{aligned}
\end{equation}
For simplicity let us break the right-hand side of equation $(\ref{equation 3.1.4})$ into two parts :
 \begin{equation}
     \Delta_{1_1} =\sum_{i=1}^q \int_{G_i}\sup_{y_1,y_2 \in S_r(x)} \left| f(\tau_i^{-1}(y_1))\bigl(g(\tau_i^{-1}(y_1))-g(\tau_i^{-1}(y_2))\bigl)\right|\, dm(x) .
 \end{equation}
 and 
\begin{align}
    \Delta_{1_2} &=\sum_{i=1}^q \int_{G_i}\sup_{y_1,y_2 \in S_r(x)} \left| g(\tau_i^{-1}(y_2))\bigl(f(\tau_i^{-1}(y_1))-f(\tau_i^{-1}(y_2))\bigl)\right| \,dm(x) .
\end{align}
First, we estimate $ \Delta_{1_{1}}$. For any $i=1,2,\dots,q$, the derivative $\tau'$ is $\varepsilon$-H\"older continuous on  $I_i$ and $\varepsilon=1/p$, with constant $M_i$ and   $\abs{\tau_i' (x)}\geq s_i >1$.  Let $M=\displaystyle\max_{1 \leq i\leq q} M_i$ and $s=\displaystyle\min_{1 \leq i\leq q} s_i$.
Also, \\$\max\{\abs{y_1-y_2}\mid y_1,y_2 \in S_r(x)\}=2r$, and $\displaystyle \abs{\tau_i^{-1}(y_1)-\tau_i^{-1}(y_2)}\leq \frac{1}{s} \cdot \abs{y_1-y_2}$. Then,
\begin{equation*}
\begin{aligned}
 \abs{ g(\tau_i^{-1}(y_1))-g(\tau_i^{-1}(y_2))}
   \leq \bigg|\frac{\tau'(\tau_i^{-1}(y_1))-\tau'(\tau_i^{-1}(y_2))}{\tau'(\tau_i^{-1}(y_1))\cdot\tau'(\tau_i^{-1}(y_2))}\bigg|\leq & \frac{2^{1/p} M\cdot(r/s)^{1/p}}{s}\cdot \frac{1}{\abs{\tau'(\tau_i^{-1}(y_1))}} .&\\
    \end{aligned}
    \end{equation*}
   For simplicity let $D=\frac{M \cdot r^\frac{1}{p}}{s^{1+\frac{1}{p}}}$. Hence,
   \begin{equation*}
       \begin{aligned}
           \Delta_{1_1}&\leq 2^{1/p} D\sum_{i=1}^q \int_{G_i}\sup_{y_1 \in S_r(x)} \left|f(\tau_i^{-1}(y_1))\right| \cdot \frac{1}{\abs{\tau'(\tau_i^{-1}(y_1))}}\,dm(x) .&\\
       \end{aligned}
   \end{equation*}
\begin{figure}[H]
    \centering
    \includegraphics[width=1\linewidth]{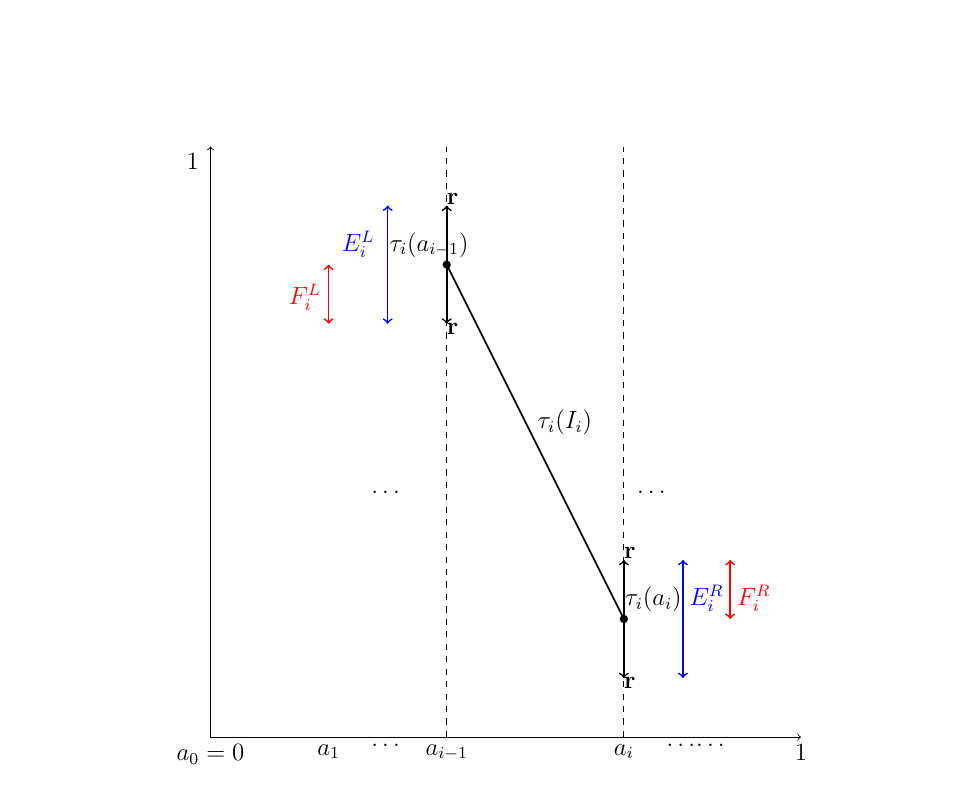}
    \caption{The sets $E_i^L, F_i^L, E_i^R$ and $F_i^R$.}
    \label{fig:enter-label}
\end{figure} 
 Adding and subtracting $f(\tau_i^{-1}(x))$ we obtain,
    \begin{align*}
        \Delta_{1_1} &\leq2^{1/p} D\sum_{i=1}^q \int_{G_i}\sup_{y_1\in S_r(x)} \left|f(\tau_i^{-1}(x))\right| \cdot \left|g(\tau_i^{-1}(y_1))\right|\,dm(x) &\\
        &+ 2^{1/p} D\sum_{i=1}^q \int_{G_i}\sup_{y_1 \in S_r(x)} \left|f(\tau_i^{-1}(y_1))-f(\tau_i^{-1}(x))\right|\cdot \left|g(\tau_i^{-1}(y_1))\right| \,dm(x)&\\
        &\leq 2^{1/p} D\sum_{i=1}^q \int_{G_i}\sup_{y_1 \in S_r(x)} \left|f(\tau_i^{-1}(x))\right| \cdot \left|g(\tau_i^{-1}(y_1))\right|\,dm(x) &\\
        &+ 2^{1/p} D\sum_{i=1}^q \int_{G_i}\sup_{y_1 \in S_r(x)}  \text{osc}(f(\tau_i^{-1}(x), r, x) \cdot \left|g(\tau_i^{-1}(y_1))\right|  \,dm(x).&
           \end{align*}
            By Proposition \ref{prop 1.7.8}, we have $ \text{osc}\left(f(\tau_i^{-1} (x)),r, x\right)\leq  \text{osc}\left(f, \frac{r}{s},\tau_i^{-1} (x)\right)$ which gives,
             \begin{align*}
          \Delta_{1_1}&\leq 2^{1/p} D\sum_{i=1}^q \int_{G_i}\sup_{y_1 \in S_r(x)} \left|f(\tau_i^{-1}(x))\right| \cdot \left|g(\tau_i^{-1}(y_1))\right|\,dm(x) &\\
        &+2^{1/p} D\sum_{i=1}^q \int_{G_i}\sup_{y_1 \in S_r(x)} \text{osc}\left(f, \frac{r}{s},\tau_i^{-1} (x)\right) \cdot \left|g(\tau_i^{-1}(y_1))\right| \,dm(x) .&
     \end{align*}
      Now, for any $w\in S_r(x)$ let us  estimate $\left|\frac{g(\tau_i^{-1}(w)) }{g(\tau_i^{-1}(x))}\right|$ :
       \begin{equation*}
         \left|\frac{g(\tau_i^{-1}(w)) }{g(\tau_i^{-1}(x))}-1\right|=\left|\frac{\frac{1}{\tau_i'(\tau_i^{-1}(w))}-\frac{1}{\tau_i'(\tau_i^{-1}(x))}}{\frac{1}{\tau_i'(\tau_i^{-1}(x))}}\right|=\left|\frac{\tau_i'(\tau_i^{-1}(x))-\tau_i'(\tau_i^{-1}(w))}{\tau_i'(\tau_i^{-1}(w))}\right|\leq \frac{M \cdot (r/s)^\frac{1}{p}}{s}.
     \end{equation*}
     So,
     \begin{equation}
        \left|\frac{g(\tau_i^{-1}(w)) }{g(\tau_i^{-1}(x))}\right|\leq 1+\frac{M \cdot (r/s)^\frac{1}{p}}{s}=1+ D.
     \end{equation}
     Then,
     \begin{equation}
          \begin{aligned}
        \Delta_{1_1} &\leq  2^{1/p} D\left(1+D\right)\sum_{i=1}^q\int_{G_i}  \left|f(\tau_i^{-1}(x))\right|\left|g(\tau_i^{-1}(x))\right|\,dm(x) &\\
        &+ 2^{1/p} D\left(1+D\right)\sum_{i=1}^q\int_{G_i} \text{osc}\left(f, \frac{r}{s},\tau_i^{-1} (x)\right)\left|g(\tau_i^{-1}(x))\right| \,dm(x) .&
     \end{aligned}
     \end{equation}
      Now, we estimate $\Delta_{1_2}$ :
        \begin{align*}
    \Delta_{1_2} &=\sum_{i=1}^q \int_{G_i} \sup_{y_1,y_2\in S_r(x)} \left| g(\tau_i^{-1}(y_2))\bigl(f(\tau_i^{-1}(y_1))-f(\tau_i^{-1}(y_2))\bigl)\right| \,dm(x)&\\
    &\leq (1+D)\sum_{i=1}^q \int_{G_i}\sup_{y_1,y_2\in S_r(x)}  \left|g(\tau_i^{-1}(x))\right|\left|f(\tau_i^{-1}(y_1))-f(\tau_i^{-1}(y_2))\right|\,dm(x) &\\
    &\leq(1+D)\sum_{i=1}^q \int_{G_i}  \text{osc}\left(f, \frac{r}{s},\tau_i^{-1} (x)\right) \left|g(\tau_i^{-1}(x))\right|\,dm(x).& 
\end{align*}
Combining the estimates for $\Delta_{1_1}$ and $\Delta_{1_2}$, we get
\begin{equation}
    \begin{aligned}
        \Delta_{1}&\leq  2^{1/p} D\left(1+D\right) \sum_{i=1}^q\int_{G_i} \left|f(\tau_i^{-1}(x))\right|\left|g(\tau_i^{-1}(x))\right|\, dm(x) &\\
        &+ 2^{1/p} D\left(1+D\right) \sum_{i=1}^q\int_{G_i} \text{osc}\left(f, \frac{r}{s},\tau_i^{-1} (x)\right)\left|g(\tau_i^{-1}(x))\right| \, dm(x) &\\
        &+(1+D) \sum_{i=1}^q\int_{G_i} \text{osc}\left(f, \frac{r}{s},\tau_i^{-1} (x)\right) \left|g(\tau_i^{-1}(x))\right|\, dm(x).
    \end{aligned}
\end{equation}
We use the change of variables, $z=\tau_i^{-1}(x)$ with $dm(z)=\frac{1}{\abs{\tau'(\tau_i^{-1}(x))}} \, dm(x) =\abs{g(\tau_i^{-1}(x))} \, dm(x)$, and obtain
    \begin{equation*}
        \begin{aligned}
          \Delta_1&\leq 2^{1/p} D\left(1+D\right) \sum_{i=1}^q\int_{I_i} \left|f(z)\right|dm(z)
        + 2^{1/p} D\left(1+D\right) \sum_{i=1}^q\int_{I_i} \text{osc}\left(f, \frac{r}{s},z\right) dm(z) &\\
        &+(1+D) \sum_{i=1}^q\int_{I_i} \text{osc}\left(f, \frac{r}{s},z\right) dm(z)&\\
        &\leq 2^{1/p} D(1+D) \norm{f}_1+ 2^{1/p} D(1+D)\,  \text{osc}_1\left(f, \frac{r}{s}\right)+(1+D) \text{osc}_1\left(f, \frac{r}{s}\right)
        \end{aligned}
    \end{equation*}
Now, let us estimate $\Delta_2$.
\begin{equation}\label{eq 3.2.10}
    \Delta_2 =
\sum_{i\in Q}
\int_{E_i^R}
\sup_{y_1,y_2\in F_i^{R\cup}}
\bigl|h_i(y_1)-h_i(y_2)\bigr|\, dm(x)
\;+\;
\sum_{i\in Q}
\int_{E_i^L}
\sup_{y_1,y_2\in F_i^{L\cup}}
\bigl|h_i(y_1)-h_i(y_2)\bigr|\, dm(x).
\end{equation}
For \(x\in E_i^R\cup E_i^L\), the supremum can be attained (or almost attained) on pairs
\(y_1,y_2\) such that both points have \(\tau_i\) preimages, or on such that one of them does
not have a \(\tau_i\) preimage. Let us consider \(x\in E_i^R\). For \(y_1,y_2\in F_i^{R\cup}\), if both have $\tau_i$ preimages, we obtain
\begin{align*}
    \left|h_i(y_1)-h_i(y_2)\right|&=\left|f\cdot g(\tau^{-1} y_1)-f\cdot g(\tau^{-1}y_2)\right|&\\
    &\leq \left|f(z_1)-f(z_2)\right| g(z_1)+ \left|f(z_2)\right| \left(g(z_1)- g(z_2)\right)&\\
    &\leq \frac{1}{s}\,\bigl|f(z_1)-f(z_2)\bigr|+ \left|f(z_2)\right| \frac{M}{s^2}\left(\frac{2r}{s}\right)^{1/p}
\end{align*}
where \(z_1=\tau_i^{-1}(y_1)\), \(z_2=\tau_i^{-1}(y_2)\in \tau_i^{-1}\!\bigl(F_i^{R\cup}\bigr)\). We used the Hölder condition for $g = \frac{1}{\tau'}$. For $y\in F_i^{R\cup}$,  we have
\begin{equation}\label{eq:hi-single}
\bigl|h_i(y)\bigr|\;\le\;\frac{1}{s}\,\bigl|f(z)\bigr|,
\end{equation}
where \(z=\tau_i^{-1}(y)\in \tau_i^{-1}\!\bigl(F_i^{R\cup}\bigr)\).
We have the corresponding estimates for \(x\in E_i^L\).

For sufficiently small \(A\) (in the definition of the oscillation), all the sets
\(\tau_i^{-1}\!\bigl(F_i^{R\cup}\bigr)\), \(\tau_i^{-1}\!\bigl(F_i^{L\cup}\bigr)\), \(i\in Q\), are disjoint, with the Lebesgue measure of each of them being less than $2r/s$. We choose disjoint sets $Y_{i}^R,Y_i^L$, $i=1,2,...,q$ such that $Y_{i}^R \supset\tau_{i}^{-1}\left(F_{i}^{R\cup}\right)$, $Y_i^L \supset\tau_i^{-1}\left(F_{i}^{L\cup}\right) $ for all $r\leq A$. We choose $A$ such that these sets are disjoint. We can choose them of equal length, say $B$, independent of $r$ and with $B=2A$. \\
   Invoking Propositions~\ref{prop 1.7.9} with $k=2$ and~\ref{Prop 1.7.10}, we see that in both cases
(two preimages or one) we have
\[
\sup_{y_1,y_2\in F_i^{R\cup}} |h_i(y_1)-h_i(y_2)|
\;\le\;\left(\frac{1}{s}+\frac{M}{s^2}\left(\frac{2r}{s}\right)^{1/p}\right)\left(
\frac{1}{A}\int_{Y_i^R} \operatorname{osc}(f,A,x)\,dm(x)
\;+\;
\frac{1}{A}\int_{Y_i^R} |f(x)|\,dm(x)\right),
\]
with corresponding estimates for $\sup_{y_1,y_2\in F_i^{L\cup}}\left|h_i(y_1)-h_i(y_2)\right|$.

Using this in \eqref{eq 3.2.10} we obtain
\begin{equation}
    \begin{aligned}
        \Delta_2 \;&\le\;
\frac{2r}{s}\left(1+2^{1/p} D\right)\sum_{i\in Q}\left(
\int_{Y_i^R\cup Y_i^L} \frac{\lvert \operatorname{osc}(f,A,z)\rvert}{A}\,dm(z)
\;+\;
\frac{1}{A}
\int_{Y_i^R\cup Y_i^L} |f(z)|\,dm(z)\right)&\\
&\leq \frac{2r}{s}\left(1+2^{1/p} D\right) \frac{\operatorname{osc}_1(f,A)}{A}+ \frac{2r}{s}\left(1+2^{1/p} D\right) \frac{1}{A} \norm{f}_1,
    \end{aligned}
\end{equation}
where $D=\frac{M r^{1/p}}{s^{1+1/p}}$, as introduced before.

Combining estimates for $\Delta_1$ and $\Delta_2$ we obtain 
      \begin{equation*}
        \begin{aligned}
             \text{osc}\,_1(P_\tau f,r) &\leq 2^{1/p} D(1+D) \norm{f}_1+ 2^{1/p} D(1+D)\,  \text{osc}_1\left(f, \frac{r}{s}\right)+(1+D) \text{osc}_1\left(f, \frac{r}{s}\right)&\\
             &+  \frac{2r}{s}\left(1+2^{1/p} D\right)\left(\frac{\text{osc}\,_1(f,A)}{A} + \frac{1}{A} \norm{f}_1\right).&
        \end{aligned}
    \end{equation*}
    Note that $\frac{D}{r^{1/p}}=\frac{M}{s^{1+1/p}}$ is independent of $r$ and $D\leq \frac{M \cdot A^{1/p}}{s^{1+1/p}}$. This allows us to make $D$ arbitrarily small by choosing appropriate $A$. Dividing by $r^{1/p}$ and taking supremum over $0<r\leq A$, we get
\begin{equation*}
    \begin{aligned}
       \text{var}_{1,\frac{1}{p}}(P_\tau f)&\leq \frac{2^{1/p} M(1+D)}{ s^{1+1/p}}\norm{f}_1+\frac{2^{1/p} D(1+D)}{s^{1/p}} \text{var}_{1,\frac{1}{p}} f+ \frac{(1+D)}{s^{1/p}} \text{var}_{1,\frac{1}{p}}f &\\
          &+\frac{2}{s}\cdot \frac{A^{1-1/p}}{A^{1-1/p}} \left(1+2^{1/p} D\right) \text{var}_{1,\frac{1}{p}} f
        +\frac{2}{s}\cdot\frac{A^{1-1/p}}{A}\left(1+2^{1/p} D\right) \norm{f}_1,&\\
        &=\left(\frac{2^{1/p} D(1+D)}{s^{1/p}}+ \frac{(1+D)}{s^{1/p}}+\frac{2\left(1+2^{1/p} D\right) }{s}\right) \text{var}_{1,\frac{1}{p}}f&\\
        &+\left(\frac{2^{1/p} M(1+D)}{ s^{1+1/p}}+\frac{2}{s}\cdot\left(1+2^{1/p} D\right) \frac{A^{1-1/p}}{A}\right)\norm{f}_1&\\
        &=\alpha \cdot \text{var}_{1,\frac{1}{p}}f+ \beta\cdot \norm{f}_1,
    \end{aligned}
\end{equation*}
where $\alpha=\frac{2^{1/p} D(1+D)}{s^{1/p}}+ \frac{(1+D)}{s^{1/p}}+\frac{2\left(1+2^{1/p} D\right) }{s}$ and $\beta=\frac{2^{1/p} M(1+D)}{ s^{1+1/p}}+\frac{2}{s}\cdot\left(1+2^{1/p} D\right) \frac{1}{A^{1/p}}$.
\end{proof}
\textnormal{Since $D$ can be made arbitrarily small, to make $\alpha<1$ we need}
\begin{equation}\label{eq 3.2.16}
    \frac{1}{s^{1/p}}+\frac{2}{s}<1.
\end{equation}
\begin{theorem}\label{thm 3.2.2}
         Let \( \tau \in \mathcal{T}_{1+\varepsilon} \) for some \( p \geq 1 \) and $\varepsilon=1/p$. If \( s=\min_{1 \leq i \leq q} s_i\) satisfies condition (\ref{eq 3.2.16}),  then we obtain the Lasota–Yorke inequality: There exist constants \( 0 < \alpha < 1 \) and \( K > 0 \) such that
          \begin{equation}
             \|P_\tau f \|_{1,\frac{1}{p}} \leq \alpha \cdot \| f \|_{1,\frac{1}{p}} + K \cdot \| f \|_1 , 
          \end{equation}
for all \( f \in BV_{1,\frac{1}{p}} \).
      \end{theorem}   
  \begin{proof}
 Recall from the proof of Lemma \ref{lem 3.2.1} that
\[
D \leq M A^{1/p} \cdot s^{-(1+1/p)}.
\]
By choosing \( A \) sufficiently small, we can ensure that \( \alpha < 1 \). For every \( f \in BV_{1,\frac{1}{p}} \), by Definition 1.3, we have
\[
\| P_\tau f \|_{1,\frac{1}{p}} = \operatorname{var}_{1,\frac{1}{p}}(P_\tau f) + \| P_\tau f \|_1.
\]
Using the properties established earlier in Lemma \ref{lem 3.2.1}, we bound the variation term:
\begin{align*}
  \operatorname{var}_{1,\frac{1}{p}} (P_\tau f) &\leq \alpha \cdot \operatorname{var}_{1,\frac{1}{p}} f + \beta \cdot \| f \|_1&\\ 
  &\leq \alpha (\norm{f}_{1,\frac{1}{p}} - \norm{f}_1)+\beta \cdot \| f \|_1.
 \end{align*}
Adding \( \| P_\tau f \|_1 =\|f\|_1\) to both sides, we obtain:
\[
\| P_\tau f \|_{1,\frac{1}{p}} \leq \alpha \cdot \| f \|_{1,\frac{1}{p}} - \alpha \cdot \norm{f}_1+\beta \cdot \| f \|_1 + \| f \|_1.
\]
Rewriting the right-hand side, we obtain:
\[
\| P_\tau f \|_{1,\frac{1}{p}} \leq \alpha \cdot \| f \|_{1,\frac{1}{p}} + K \cdot \| f \|_1,
\]
where $K=(1-\alpha+\beta)>0$.
\end{proof}
 \begin{theorem}\label{thm 3.2.3}
      Let $f\in BV_{1,\frac{1}{p}}$, $1\le p <\infty$. There exist a constant $C>0$ such that for every sufficiently large $n\in \mathbb{N}$,
\begin{equation*}
    \norm{P^n_\tau f}_{1,\frac{1}{p}}\leq C\cdot \norm{f}_1,
\end{equation*}
  where $C=\left(1+\frac{K}{1-\alpha}\right) $, $\alpha<1$.
     \end{theorem}
     \begin{proof} From Theorem \ref{thm 3.2.2}, with the appropriate $0<\alpha<1$ and $K>0$ we have,
     \[
\| P_\tau f \|_{1,\frac{1}{p}} \leq \alpha \cdot \| f \|_{1,\frac{1}{p}} + K \cdot \| f \|_1,
\]
for all $f\in BV_{1,\frac{1}{p}}$. Replacing $f$ by $P_\tau f$ in above equation yields
 \begin{equation}\label{equ 3.2.18}
  \begin{aligned}
       \norm{P_\tau^2 f}_{1,\frac{1}{p}}&\leq \alpha \cdot \norm{P_\tau f}_{1,\frac{1}{p}}+ K \cdot \norm{P_\tau f}_1&
  \end{aligned}
 \end{equation}
since $P_\tau$ is a contraction on $L^1$, and using equation (\ref{equ 3.2.18}), we obtain,
\begin{equation}
    \begin{aligned}
          \norm{P_\tau^2 f}_{1,\frac{1}{p}}&\leq\alpha \cdot\left( \alpha \cdot \| f \|_{1,\frac{1}{p}} + K \cdot \| f \|_1 \right)+K \cdot \norm{f}_1&\\
       &=\alpha^2 \cdot \norm{f}_{1,\frac{1}{p}} +(K + K \alpha) \cdot\norm{f}_1.
    \end{aligned}
\end{equation}
 Similarly, for $P_\tau^3 f$ we get,
 \begin{equation}
      \norm{P_\tau^3 f}_{1,\frac{1}{p}}\leq \alpha^3 \cdot  \norm{f}_{1,\frac{1}{p}}+(K + K \alpha +  K \alpha^2) \cdot \norm{f}_1.
 \end{equation}
By induction, for $n\in \mathbb{N}$
 \begin{equation*}
 \begin{aligned}
     \norm{P_\tau^n f}_{1,\frac{1}{p}}&\leq \alpha^n \cdot \norm{f}_{1,\frac{1}{p}} +(K + K \alpha +  K \alpha^2+\dots + K \alpha^{n-1})\cdot \norm{f}_1&\\
     &=\alpha^n \cdot \norm{f}_{1,\frac{1}{p}}+K (1 +  \alpha +   \alpha^2+\dots +  \alpha^{n-1})\cdot \norm{f}_1 ,&\\
 \end{aligned}
 \end{equation*}
  so by induction, for sufficiently large  $n\in \mathbb{N}$,
  \begin{equation}
       \norm{P_\tau^n f}_{1,\frac{1}{p}}\leq \left(1+\frac{K}{1-\alpha}\right)\norm{f}_1, 
  \end{equation}  
  and $C=1+\frac{K}{1-\alpha}$.
     \end{proof}
     \textnormal{The property of $P_\tau$ and the space $BV_{1,\frac{1}{p}}$ we proved in Theorem \ref{thm 3.2.2} and Theorem \ref{thm 1.7.7} allow us to use  \cite{ionescu1950}, with $X=BV_{1,\frac{1}{p}}$, $Y=L^1$, and $P=P_\tau$ which in particular states that on each ergodic component $\tau$ has ACIM.}
\begin{theorem}\label{thm 3.2.4}
    Under the assumptions of Theorem \ref{thm 3.2.2}, the following hold:
    \begin{enumerate}
        \item $P_\tau: BV_{1,\frac{1}{p}} \rightarrow BV_{1,\frac{1}{p}}$ has a finite number of eigenvalues $c_1,c_2,...c_r$ of modulus $1$;
        \item $E_i=\{f\in L^1 \mid P_\tau f=c_if\} \subseteq BV_{1,\frac{1}{p}}$, and $E_i$ is finite-dimensional for $i=1,2..,r$;
        \item $\displaystyle P_\tau =\sum_{i=1}^r c_i \Psi_i +Q$, where $\Psi_i$ represents the projection onto the eigenspace $E_i$, $\norm{\Psi_i}_1= 1$, and $Q$ is a linear operator on $L^1$ with $Q(BV_{1,\frac{1}{p}}) \subseteq BV_{1,\frac{1}{p}}$, $\displaystyle \sup_{n\in \mathbf{N}} \norm{Q^n}_1 < \infty$ and $\norm{Q^n}_{1,\frac{1}{p}}= O(q^n)$ for some $0<q<1$. Furthermore, $\Psi_i \Psi_j=0 (i \neq j) $ and $\Psi_i Q= Q \Psi_j=0$ for all $i$.
        \end{enumerate}
\end{theorem}
\begin{proof}Since $P_\tau$ is a contraction on $L^1$, its powers are uniformly bounded. That is, we have $\|P_\tau^n\|_{1} \leq 1$ for all $n \in \mathbb{N}$.

From Theorem~\ref{thm 3.2.2}, we know that the operator satisfies a Lasota–Yorke inequality of the form $$\|P_\tau f\|_{1, \frac 1 p} \leq \alpha \|f\|_{1, \frac 1 p} + K \|f\|_{1},
$$
for some constant $K > 0$ and $0<\alpha<1$. In addition, by Theorem~\ref{thm 1.7.7}, the embedding $
BV_{1,\frac{1}{p}} \hookrightarrow L^1
$ is compact, hence the operator $P_\tau: (BV_{1,\frac{1}{p}}, \|\cdot\|_{{1,\frac{1}{p}}}) \to (L^1, \|\cdot\|_1)$ is compact. These properties shows that $P_\tau$ satisfies the hypotheses of the \cite{ionescu1950} theorem, which gives quasi-compactness of $P_\tau$ on $BV_{1,1/p}$ and provides the stated spectral decomposition: finitely many eigenvalues of modulus $1$ with finite-dimensional eigenspaces, and a remainder $Q$ whose iterates are uniformly bounded in $L^1$ and contract exponentially fast in the $BV_{1,1/p}$ norm. The projection properties follow from the spectral decomposition. This proves (1)--(3).
    \end{proof}
    
   \textnormal{ Note that $P_\tau$, is a positive operator from $L^1$ to $L^1$, and we can use Rota's theorem to prove the proposition below.}
    \begin{proposition}\label{prop 3.2.5}
    Let $\tau \in \mathcal T_{1+\varepsilon}$ for some $\varepsilon >0$. Then $\tau$ has a finite number of ergodic components. On each component, $\tau$ has an ACIM. On each component,  some iterate of $\tau$, say $\tau^n$, has a finite number of disjoint invariant domains and on each of them $\tau^n$ is exact.
\end{proposition}
\begin{proof}
 By the Ionescu-Tulcea and Marinescu theorem, $P_\tau$ has a finite number of eigenvalues of modulus
1. By Rota's theorem they form a group or a union of groups of roots of unity. Each such group
corresponds to an ergodic component of $\tau$. One of the eigenvalues in the group is 1, which
shows the existence of an ACIM. If the order of roots in the group is $n$, then $\tau^n$
has all eigenvalues equal to 1, and the supports of the eigenfunctions form the disjoint
$\tau^n$-invariant domains. On each domain $\tau^n$ is exact. The exactness for the case when $\tau$ has unique ACIM $\mu=h \cdot m$, follows from $3$ of Theorem \ref{thm 3.2.4}. If $f\in BV_{1,\frac{1}{p}}$ then
\begin{equation*}
    P_\tau^n f =\int_I f dm \cdot h+ Q^n f \rightarrow \int_I f dm \cdot h,
\end{equation*}
as $n\rightarrow \infty$, both in $BV_{1,\frac{1}{p}}$ and $L^1$. Since $BV_{1, \frac{1}{p}}$ is dense in $L^1$ and $P_\tau$ is a contraction on $L^1$ the convergence in $L^1$ holds for all functions in $L^1$. And this implies the exactness.
\end{proof}
\begin{proposition}\label{prop:roots-of-unity}
Let $c$ with $|c|=1$ be an eigenvalue of $P_{\tau}$ with corresponding eigenvector
$f\in BV_{1,\frac 1 p}$. Then, for every $n\in\mathbb{N}$ with $n\ge 2$, the number $c^{\,n}$ is also
an eigenvalue of $P_{\tau}:BV_{1, \frac 1 p}\to BV_{1, \frac 1 p}$.
\end{proposition}
\begin{proof}
By the argument in the proof of Rota’s theorem in \cite{schaefer1964} , if $f=|f|\cdot g$ (so $|g|=1$ a.e.), then
for every $n\ge 1$,
\[
c^{\,n} f \;=\; P_{\tau}\!\big(|f|\cdot g^{\,n}\big).
\]
Thus it is enough to show that $|f|\cdot g^{\,n}\in BV_{1, \frac 1 p}$.
Since $||f|(x)-|f|(y)|\le |f(x)-f(y)|$ for all $x,y$, we have $|f|\in BV_{1,\frac 1 p}$.
From the same proof of Rota’s theorem, $|f|$ is the density of an absolutely continuous
invariant measure for $\tau$. Because $\tau$ is piecewise expanding with bounded derivative,
a standard result (see, e.g., \cite{keller1985}) implies that $|f|$ is bounded away from $0$.
Hence $g=\frac {f} {|f|}\in BV_{1, \frac 1 p}$, and the product rule for $BV_{1, \frac 1 p}$ with bounded multipliers
gives $|f|\cdot g^{\,n}\in BV_{1,\frac 1 p}$. The claim follows.
\end{proof}
\begin{corollary}\label{cor:roots-union}
The eigenvalues $c_1,c_2,\dots,c_r$ listed in Theorem~\ref{thm 3.2.4}
form a finite union of groups of roots of unity.
\end{corollary}
\begin{proof}
Combine Proposition~\ref{prop:roots-of-unity} with statement 1. of Theorem~ \ref{thm 3.2.4}.
\end{proof}
    \begin{proposition}\label{Prop 3.2.9}
    Let \( \tau \in \mathcal{T}_{1+\frac 1 p} \) for some \( p \geq 1 \). Then $\tau$ has a finite number of ergodic components. On each component, $\tau$ has an ACIM. On each component,  some iterate of $\tau$, say $\tau^n$, has a finite number of disjoint invariant domains and on each of them $\tau^n$ is exact.
\end{proposition}
\begin{proof}
 By the Ionescu-Tulcea and Marinescu theorem, $P_\tau$ has a finite number of eigenvalues of modulus
1. By Rota's theorem they form a group or a union of groups of roots of unity. Each such group
corresponds to an ergodic component of $\tau$. One of the eigenvalues in the group is 1, which
shows the existence of an ACIM. If the order of roots in the group is $n$, then $\tau^n$
has all eigenvalues equal to 1, and the supports of the eigenfunctions form the disjoint
$\tau^n$-invariant domains. On each domain $\tau^n$ is exact. The exactness for the case when $\tau$ has unique ACIM $\mu=h \cdot m$, follows from $(3)$ of Theorem \ref{thm 3.2.4}: If $f\in BV_{1,\frac{1}{p}}$ then
\begin{equation*}
    P_\tau^n f =\int_I f dm \cdot h+ Q^n f \rightarrow \int_I f dm \cdot h,
\end{equation*}
as $n\rightarrow \infty$, both in $BV_{1,\frac{1}{p}}$ and $L^1$. Since $BV_{1, \frac{1}{p}}$ is dense in $L^1$ and $P_\tau$ is a contraction on $L^1$ the convergence in $L^1$ holds for all functions in $L^1$. And this implies the exactness.
\end{proof}
\begin{proposition}\label{prop 3.2.6}
Let \( \tau \in \mathcal{T}_{1+\frac 1 p} \) for some \( p \geq 1 \). Then $P_\tau$, considered as an operator from $L^1$ to $L^1$, does not have any eigenvalues of modulus 1 except those listed in Theorem \ref{thm 3.2.4}.
\end{proposition}
\begin{proof}
     All eigenfunctions listed in  Theorem \ref{thm 3.2.4} are in $BV_{1, \frac{1}{p}}$. By the results of \cite{schaefer1964}, if $f\in L^1$ is an eigenfunction of $P_\tau$, then $|f|$ is a density of an acim. If $f\in BV_{1, \frac{1}{p}}$, then $|f|\in BV_{1, \frac{1}{p}}$, as well.  Thus, by \cite{keller1978} their supports are finite unions of closed intervals.
If $S$ is the union of their supports, then $S$ is also a finite union of closed intervals. 

Let $g\in L^1$ be an eigenfunction of $P_\tau$ not listed in Theorem \ref{thm 3.2.4}.
First, we will prove that $S_g=\supp g \cap (I\setminus S)$ is of measure 0.
Let us assume that $m(S_g)>0$ and let $J\supset S_g$ be an interval contained in $I\setminus S$.
 Let $f=\chi_J$. Then, by Theorem \ref{thm 3.2.3},  the iterates $P^n_\tau f $ are  uniformly bounded in $BV_{1, \frac{1}{p}}$. By Yosida-Kakutani Theorem \cite{yosidakakutani1941},
$\tau$ preserves an ACIM with support containing the set $S_g$. By Theorem~\ref{thm 1.7.7}, its density is in $BV_{1, \frac{1}{p}}$. This contradicts the definition of $S$. 

We proved that $S_g\subset S$. Let $k\in\mathbb N$ be such that all the eigenvalues listed in Theorem \ref{thm 3.2.4} satisfy $c_i^k=1$, $k=1,2,\dots,r$. Let $J_1,\dots, J_p$ be the ergodic components for $\tau^k$. Then, $\tau^k$ is exact on each set $J_i$ and  $S=\bigcup_{i=1}^p J_i$. 
Let $f_i\in BV_{1, \frac 1p}$ be the density of the acim on $J_i$, $i=1,2,\dots,p$.
Let $P_\tau g=\lambda g$, i.e., the eigenfunction $g$ corresponds to the eigenvalue $\lambda$.
 Let $i$ be one of the numbers $\{1,2,\dots,p\}$. Since $\tau^k$ is
 exact on $J_i$, we know that for any $f\in L^1$, $P_{\tau^k} ^{n} f\to \int_{J_i} f dm \cdot f_i$, as $n\to \infty$.
In particular, $P_{\tau^k} ^{n} g_i =(\lambda^k)^n g_i=g_i \to \int_{J_i} g_i dm \cdot f_i$,
and $g_i$ is a multiple of $f_i$, and $g_i$ is already listed.
This completes the proof.
\end{proof}
\textnormal{The following two theorems assume that $\tau\in \mathcal{T}_{1+\varepsilon}$ for some $\varepsilon=\frac{1}{p}$ and \( p \geq 1 \), admits a unique ACIM.}
\begin{theorem}\label{Thm 3.2.7}
    Let $\tau \in \mathcal T_{1+\varepsilon}$ for some $\varepsilon >0$,  and let $\mu=h \cdot m$ be the unique absolutely continuous $\tau$-invariant measure. For $f \in BV_{1,\frac{1}{p}}$ and $g \in L^\infty$, the correlation $C_m(f, g, N)$ decays exponentially as $N \rightarrow \infty$.
\end{theorem}
\begin{proof}
    There is only one eigenvalue of modulus $1$, which is $1$. For $f \in BV_{1,\frac{1}{p}}$ and $g \in L^\infty$, the correlation after $N$ iterations is given by
    \begin{equation*}
        C_m(f, g, N) = \left| \int_I f \cdot g \circ \tau^N dm-\int_I f dm \int_I g d\mu \right|=\bigg|\int_I P_{\tau}^N f \cdot g \, dm - \int_I f \, dm \int_I g \, d\mu\bigg|.
    \end{equation*}
    By the Ionescu Tulcea and Marinescu theorem, we know that
    \begin{equation*}
        P_{\tau}f = \Psi_1 f + Q f,
    \end{equation*}
    where $\Psi_1f$ is the projection of $f$ on the eigenfunction which is the invariant density $h$ so the projection $\Psi_1f=\int_I f dm \cdot h=m(f)\cdot h$ and $\norm{Q^N}_{1,\frac{1}{p}} = O(q^N), 0<q< 1$. Then, for any $N$,
    \begin{equation*}
        P_{\tau}^Nf = \Psi_1f + Q^Nf.
    \end{equation*}
    Hence, the correlation becomes
    \begin{equation*}
        \begin{aligned}
            C_m(f, g, N) &= \bigg| \int_I \Psi_1f \cdot g \, dm + \int_I Q^Nf \cdot g \, dm - m(f) \cdot \mu(g)\bigg| \\
            &= \bigg|\int_I m(f) \cdot h \cdot g \, dm + \int_I Q^Nf \cdot g \, dm - m(f) \cdot \mu(g)\bigg| \\
            &= \bigg|m(f) \mu(g) + \int_I Q^Nf\cdot g \, dm - m(f) \mu(g)\bigg| \\
            &= \bigg|\int_I Q^Nf \cdot g \, dm\bigg|.
        \end{aligned}
    \end{equation*}
    Since $\norm{Q^N f}_{1,\frac{1}{p}} \leq C \cdot q^N \norm{f}_{1,\frac{1}{p}}$, we have
    \begin{equation*}
        C_m(f, g, N) \leq C \cdot q^N \cdot\norm{f}_{1,\frac{1}{p}} \norm{g}_\infty.
    \end{equation*}
    Since $q < 1$, $C_m(f, g, N)$ decays exponentially as $N \to \infty$.
    \end{proof}
   \textnormal{ The system $(\tau,\mu)$ has a unique absolutely continuous invariant measure $\mu$. The correlation $C_m(f, g)$ is more symmetric when the system is considered with respect to the invariant measure.}
\begin{theorem}\label{thm 3.2.8}
    Let $\tau \in \mathcal T_{1+\varepsilon}$ for some $\varepsilon >0$,  and let $\mu$ be the unique absolutely continuous $\tau$-invariant measure. For $f \in BV_{1,\frac{1}{p}}$ and $g \in L^\infty$, the correlation $C_\mu(f, g, N)$ decays exponentially as $N \rightarrow \infty$.
\end{theorem}
\begin{proof}
    The proof is similar to the proof of the previous theorem.
    \end{proof}
    \section{Frobenius-Perron operator on the space \texorpdfstring{$BV_{t,\frac{1}{p}}$}{BV(t,1/p)}}

    \textnormal{In this section, we  collect some facts about $P_\tau$ acting on the space $BV_{t, \frac{1}{p}}$, for $1\le t \le p$. Perhaps, in the future, they can be used to establish the quasi-compactness of $P_\tau$ on this spaces. Let $\mathcal{P}=\{I_i=(a_{i-1},a_i),i=1,2,3...q\}$ be a  finite partition of $I=[0,1]$ and $\varepsilon>0$. }
      \begin{lemma}\label{lem 3.3.1}
    Let $\tau \in \mathcal{T}_{1+\varepsilon}$, where $\varepsilon=1/p$ and let $1<t\le p$. Then there exist $\alpha, \beta >0$ such that, for every $f\in BV_{t,\frac{1}{p}}$,
    \begin{equation}
       \textnormal{var}_{t,\frac{1}{p}}(P_\tau f) \leq \alpha \cdot\norm{f}_{t,\frac{1}{p}} + \beta \cdot \norm{f}_t.
    \end{equation}
\end{lemma}
\begin{proof}
The proof is very similar to that of Lemma \ref{lem 3.2.1}, except for the last part where we use Proposition \ref{Prop 1.7.10} instead of \ref{prop 1.7.9}. In the end, we estimate the $L^t$ norms on each subinterval by the $L^t$ norm on the whole interval $I$. This results in the factors $q$ (the number of branches) in the estimates, and we obtain $\alpha=2q\left(\frac{2^{1/p} D(1+D)}{s^{1/p}}+ \frac{(1+D)}{s^{1/p}}+\frac{2\left(1+2^{1/p} D\right) }{s}\right)$ and $\beta=2q\left(\frac{2^{1/p} M(1+D)}{ s^{1+1/p}}+\frac{2}{s}\cdot\left(1+2^{1/p} D\right) \frac{1}{A^{1/p}}\right)$. The standard technique for lowering  $\alpha$ is to increase $s$ by using an iterate of $\tau$. This does not work when we have the number of branches in the
estimate, since then the number of branches also grows, usually exponentially.
\end{proof}
\begin{proposition}
Let \(\mu = f \cdot \nu\). Then the Frobenius-Perron operator with respect to \(\mu\) satisfies
\[
P_{\tau,\mu}g = \frac{P_{\tau,\nu}(f \cdot g)}{f}.
\]
\end{proposition}

\begin{corollary}
If \(\mu\) is \(\tau\)-invariant, then for any non-negative \(f \in L^p\), with \(1 \le p \le +\infty\), we have
\[
\int (P_{\tau,\mu} f)^p \, d\mu \le \int f^p \, d\mu.
\]
\end{corollary}

\begin{proposition}\label{prop 1.6.13}
    If $\mu=f\cdot \nu$ is $\tau$-invariant and for constants $M_1$, $M_2$ we have
 $$ f\le M_1, \ \frac 1f \le M_2 , $$
then for any non-negative $g\in L^p$, $1\le p< +\infty$, we have
$$\int (P_{\tau,\nu} g )^p d\nu \le M_1^p M_2^p \int g^p d\nu. $$
The same inequality holds for any power $P_{\tau,\nu}^k$ of the Frobenius-Perron operator.
\end{proposition}
\begin{proof}
We have
\begin{equation*}\begin{split}
\int (P_{\tau,\nu} g )^p d\nu &=\int\left(P_{\tau,\mu}\left(\frac 1f\cdot g\right) f\right)^p d\nu
=\int\left(P_{\tau,\mu}\left(\frac 1f\cdot g\right) \right)^p f^p\frac 1f d\mu\\
&\le \int\left(M_2 P_{\tau,\mu} g \right)^p M_1^{p-1} d\mu\le M_2^pM_1^{p-1}\int g^p d\mu\\
&=M_2^pM_1^{p-1}\int g^p f d\nu\le M_2^pM_1^{p}\int g^p d\nu.
\end{split}
\end{equation*}
\end{proof}
 \textnormal{ We proved in Proposition \ref{prop 3.2.5} that on any ergodic component $\tau$ has an ACIM with a density $h\in BV_{1,\frac{1}{p}}$. This implies that $h$ is bounded and has a lower-semi continuous version. We also assume that $\tau$ is piecewise expanding  with $\tau'$ bounded, as piecewise Hölder continuous function. These four facts imply, that $h$ is also bounded away from $0$. The proof is given in  \cite{keller1978}. This allows us to use Proposition \ref{prop 1.6.13} and obtain :}\\
 
  \begin{proposition}\label{prop 3.3.2}
      For some constant $L$, any $f\in BV_{t,\frac{1}{p}}$ and any $n\geq 1$, we have
      \begin{equation}
          \norm{P_\tau^n f}_t\leq L \cdot \norm{f}_t ,
      \end{equation}
      i.e., all powers of $P_\tau$ are equicontinuous.
  \end{proposition}
  \begin{theorem}\label{thm 3.3.3}
          Let \( \tau \in \mathcal{T}_{1+\varepsilon} \) for some \( 1<t\le p\) and $\varepsilon=1/p$. There exist constants \(\alpha >0 \) and \( K > 0 \) such that
          \begin{equation}
             \|P_\tau f \|_{t,\frac{1}{p}} \leq \alpha \cdot \| f \|_{t,\frac{1}{p}} + K \cdot \| f \|_t , 
          \end{equation}
for all \( f \in BV_{t,\frac{1}{p}} \). Thus, $P_\tau (BV_{t,\frac{1}{p}}) \subset BV_{t,\frac{1}{p}}$.
      \end{theorem}   
  \begin{proof}
 For every \( f \in BV_{t,\frac{1}{p}} \), by Definition \ref{def 1.7.3}, we have
\[
\| P_\tau f \|_{t,\frac{1}{p}} = \operatorname{var}_{t,\frac{1}{p}}(P_\tau f) + \| P_\tau f \|_t.
\]
Using the properties established earlier in Lemma \ref{lem 3.3.1}, we bound the variation term:
\begin{align*}
  \operatorname{var}_{t,\frac{1}{p}} (P_\tau f) &\leq \alpha \cdot \operatorname{var}_{t,\frac{1}{p}} f + \beta \cdot \| f \|_t&\\ 
  &\leq \alpha (\norm{f}_{t,\frac{1}{p}} - \norm{f}_t)+\beta \cdot \| f \|_t.
 \end{align*}
Adding \( \| P_\tau f \|_t \) to both sides, we obtain:
\begin{align*}
    P_\tau f \|_{t,\frac{1}{p}} &\leq \alpha \cdot \| f \|_{t,\frac{1}{p}} - \alpha \cdot \norm{f}_t+\beta \cdot \| f \|_t +\| P_\tau f \|_t &\\
    &\leq \alpha \cdot \| f \|_{t,\frac{1}{p}} - \alpha \cdot \norm{f}_t+\beta \cdot \| f \|_t + L \cdot \| f \|_t .
\end{align*}
Rewriting the right-hand side, we obtain:
\begin{align*}
   P_\tau f \|_{t,\frac{1}{p}} &\leq \alpha \cdot \| f \|_{t,\frac{1}{p}} + \left(L-\alpha+\beta\right) \cdot \| f \|_t&\\
   &\leq \alpha \cdot \| f \|_{t,\frac{1}{p}} + \left(1-\alpha+\beta+L\right) \cdot \| f \|_t,
\end{align*}
where $K=(1-\alpha+\beta+L)>0$.
\end{proof}
\textbf{Acknowledgments:} The research for this publication was supported by NSERC (National Science and Engineering Research Council of Canada) grant,
number RGPIN-2020-06788. The authors are very grateful to Dr. A. Berger and an anonymous reviewer for comments on the previous versions of the paper.

\end{document}